\documentclass{article}
\usepackage{graphicx} 
\newtheorem{defi}{\textbf{Definition}}[section]
\newtheorem{thm}[defi]{\textbf{Theorem}}

\newtheorem{prop}[defi]{\textbf{Proposition}}
\newtheorem{cor}[defi]{\textbf{Corollary}}
\newtheorem{sthm}{Theorem}

\newenvironment{proof}[1][Proof]{\par\noindent\textbf{#1.} }{\hfill$\square$\par}
\usepackage{amssymb,graphicx,epsfig,psfrag,epic,eepic,latexsym}
\usepackage{amsmath}
\usepackage{hyperref}
\usepackage{mathrsfs}
\usepackage{tikz}
\usepackage{amsmath}
\usepackage{comment}
\DeclareMathOperator{\arccot}{arccot}

\title{Coordinate Ordering and Non-collinear Central Configurations}
\author{Jiexin Sun\textsuperscript{}\thanks{E-mail: jiexinsun@dlut.edu.cn}
\quad Lei Zhao\textsuperscript{}\thanks{E-mail: zhao1899@dlut.edu.cn}
\\ \\
\noindent
\textsuperscript{}School of Mathematical Sciences, Dalian University of Technology,\\
Dalian, Liaoning 116024,
The People's Republic of China
}
\date{}

\begin{document}

\maketitle

\begin{abstract}
We establish a mass-independent criterion for N-body central configurations based on a special coordinate ordering of the bodies. Moreover, we identify a region in the plane that contains no non-collinear central configurations. This provides a way to exclude non-collinear central configurations that are nearly collinear.  
    
\end{abstract}

\section{Introduction}
Consider $N$ particles with positive masses $m_1,\dots,m_N$ in $\mathbb{R}^d$, $d\leq N-1$. We denote the i-th particle by $q_i=(x_i^1,\dots ,x_i^d)^T$ and the configuration by $q=(q_1,\dots,q_N)^T$. Denote $x^s=(x_1^s,\dots,x_N^s)^T,\ s=1,\dots,d$ the s-th coordinate vector of the configuration $q$. This way $q$ can be written as $q=(x^1,\dots,x^d)$. In our study, we fix the center of mass at the origin, i.e., $\sum_{i=1}^Nm_iq_i=0$.

Let 
$$
U(q)=\sum_{1\leq i<j\leq N}\frac{m_im_j}{\|q_i-q_j\|}
$$
be the force function of $N$ particles. A configuration is called a \emph{central configuration} (CC) of the Newtonian $N$-body problem if the acceleration of each particle is proportional to its position vector, which means that there exists a constant $\lambda>0$ such that
\begin{equation}\label{cceq}
  \sum_{j \neq i} \frac{m_j(q_j-q_i)}{\|q_i-q_j\|^3}
  = -\lambda\, q_i, \qquad i=1,\ldots,N.
\end{equation}
The central configuration equation (\ref{cceq}) can be written in the matrix form
\begin{equation}
Z(x^1,\dots, x^d)=\lambda (x^1,\dots, x^d),
\end{equation}
where $Z$ is a matrix with elements
$$
Z_{ij}=\left\{
\begin{aligned}
     &-\frac{m_j}{\|q_j-q_i\|^3},\quad j\neq i;  \\
     &\sum_{j\neq i} \frac{m_j}{\|q_j-q_i\|^3},\quad j=i.
\end{aligned}
\right.
$$

In other words, the $d$ coordinate vectors of a central configuration are the eigenvectors associated with the eigenvalue $\lambda$ of $Z$. The \textit{Brehm-Wintner-Conley} (BWC) matrix $Z$ depends only on the mutual distances of the particles of the configuration $q$. The matrix $Z$ is rotation invariant.

For any three positive masses given, in 1767, Euler \cite{Eu} proved the existence and uniqueness of collinear central configuration up to the order of the bodies. Lagrange \cite{Lag} proved that the only non-collinear 3-body central configurations are equilateral triangles in 1772. In 1910, Moulton \cite{Moul} extended Euler's result to any $N\geq 4$ given positive masses on a line. In \cite{Pac}, Pacella used the idea from Conley to show that we can not obtain non-collinear central configurations from a collinear central configuration under a small perturbation. 

It is natural to ask whether there are central configurations which are not collinear but `` almost collinear ". There are some mass independent criteria of not being a central configuartion. 

The $45^{\circ}$ Theorem is proved by Moeckel \cite{Moe}. He proves that there are no non-collinear central configurations with $0<\theta (q)<45^{\circ}$, where 
\begin{equation*}
\theta(q)=\min_{L}{\theta_L(q)},\quad \theta_L(q)=\max_{i\neq j}\angle(q_i-q_j,L).
\end{equation*}
Here $\angle(q_i-q_j,L)$ denotes the angle between $L$ and the line segment connecting $q_i$ and $q_j$.

The Perpendicular Bisector Theorem is a result from Conley \cite{Moe}\cite{Pac}. For any two particles $q_i$ and $q_j$ of a central configuration, the line through $q_i$ and $q_j$, together with the perpendicular bisector of the segment $q_iq_j$, divides the plane into four open quadrants. Then the remaining particles cannot all lie in the same pair of diagonally opposite quadrants.

The Disk Sector Theorem is proved by Albouy and the first-named author in \cite{AS}. For any particle $q_i$ of a planar central configuration with multiplier $\lambda$ and total mass $\mathbb{M}$, an arbitrary line passing $q_i$ and a circle with radius $\frac{\lambda}{\mathbb{M}}$ centered at $q_i$ divide the plane into four open regions. The remaining particles cannot all lie in the union of the region inside the circle on one side of the line and the region outside the circle on the other side. 

In this note, we established a new criterion for central configurations based on the ordering of coordinates. We prove
\begin{sthm}\label{ta}
    If a planar $N$-body central configuration can be written in ordered $x$ and $y$ coordinates
$$
x_1<x_2<\cdots<x_N \quad \text{and}\quad y_1<y_2<\cdots<y_N
$$
in some reference frame $xOy$, then it is collinear.
\end{sthm}

This criterion is generalized to high dimension in Proposition \ref{hdim}. A precise neighborhood of collinear configurations containing no non-collinear central configurations is obtained by this argument in Corollary \ref{MinMax}. These will be proved in Section 2. In Section 3, we compare these mass-independent central configuration criteria and propose some further questions.


\section{Main Results}

Denote $E=\{v\in \mathbb{R}^N | \sum_{i=1}^{N} m_iv_i=0 \}$. We define
$$V_+=\{v \in E: v_{1}< v_{2}<\cdots < v_{N}\},$$
$$V_-=\{v \in E: v_{1}> v_{2}>\cdots > v_{N}\},$$
$$
V_{i}=\{v_{1} < \cdots < \cdots v_{i}=v_{i+1} < \cdots < v_{N}\},\quad i=1,\dots,N-1.
$$
Then the boundary
$
\partial V_{+}=\overline{\cup V_{i}}$.

In a planar coordinate system $xOy$, we denote an $N$-body configuration by $q=(x,y)$, where $x=(x_1,\dots,x_N)^T$ and $y=(y_1,\dots,y_N)^T$. With these notations, Theorem \ref{ta} can be described more precisely as the following theorem.

\begin{thm}\label{main}(Coordinate ordering criterion)
Let $q=(x,y)$ be a central configuration.
If $x,y \in V_{+}\cup V_-$, then $q$ must be collinear.
\end{thm}

\begin{proof}
    By reflection, it suffices to prove the case $x,y\in V_+$. We consider $U=U(q)$ and $Z=Z(q)$. They only depend on the configuration $q=(x,y)$. Since $q$ is central, for some $\lambda>0$, we have
\begin{equation} \label{ccxy}
Zx=\lambda x,\quad Zy=\lambda y.
\end{equation}

We consider the gradient flow of $\mathcal{U}(v)=\frac{1}{2}\langle Mv,Zv\rangle=\frac{1}{2}v^TZ^TMv$,
\begin{equation} \label{fleq}
\dot{v}=Z v.
\end{equation}
In coordinates,
$$
\dot{v}_i=\sum_{k\neq i}\frac{m_k}{\| q_i-q_k\|^3}(v_i-v_k),\quad i=1,\dots,N.
$$
Thus
$$
\frac{d}{dt}(v_i-v_{i+1})=\sum_{k\neq i,i+1}m_k \left( \frac{v_i-v_k}{\|q_i-q_k\|^3}-\frac{v_{i+1}-v_k}{\|q_{i+1}-q_k\|^3} \right) +\frac{m_i+m_{i+1}}{\|q_i-q_{i+1}\|^3}(v_i-v_{i+1}).
$$
In particular, for $\forall v\in V_i,\ i=1,\dots,N-1$, we have
$$
\frac{d}{dt}(v_i-v_{i+1})=\sum_{k\neq i,i+1}m_k(v_i-v_k) \left( \frac{1}{\|q_i-q_k\|^3}-\frac{1}{\|q_{i+1}-q_k\|^3} \right).
$$
Since $v\in V_i$ and $x,y\in V_+$, we obtain that

\begin{itemize}
    \item for $k<i$, $\|q_i-q_k\|<\|q_{i+1}-q_k\|$;
    \item for $k>i$, $\|q_i-q_k\|>\|q_{i+1}-q_k\|$.
\end{itemize}

Therefore for $v\in V_i,\ i=1,\dots,N-1$, we have
\begin{equation} \label{vi}
\frac{d}{dt}(v_i-v_{i+1})>0.
\end{equation}
Geometrically, the gradient flow on $\partial V_{+}$ points outward from $V_+$.

Denote by $E_{xy} \subset E$ the $\mathbb{R}$-linear subspace spanned by $x$ and $y$. If $q$ is collinear, there exists $r\in \mathbb{R}\backslash \{0\}$ such that $y=rx$, so $E_{xy}\subset (V_+\cup V_-\cup \{0\})$ and $E_{xy}\cap \partial V_+=\{0\}$. 

When $q$ is non-collinear, then $E_{xy}\cap \partial V_+$ is nontrivial and there exists nonzero $\tilde{v}\in E_{xy}\cap(\cup V_i)$. Since $x,y$ satisfy (\ref{ccxy}), so is $\tilde{v}$. Together with (\ref{fleq}), $\tilde{v}$ satisfies
$$
Z\tilde{v}=\lambda \tilde{v}, \quad \dot{\tilde{v}}=\lambda \tilde{v}.
$$
Hence we have
\begin{equation}\label{c}
\frac{d}{dt}(\tilde{v}_i-\tilde{v}_{i+1})=\lambda (\tilde{v}_i-\tilde{v}_{i+1})=0,
\end{equation}
which contradicts (\ref{vi}). Therefore, $q$ has to be collinear.

\end{proof}

This argument can be easily generalized to high dimension with the same proof.

\begin{prop}\label{hdim}
    For a central configuration $q=(x^1,\dots, x^d)$ in $\mathbb{R}^d$, if $x^1,\dots, x^d\in V_+\cup V_-$, then $q$ is collinear.
\end{prop}

For a planar configuration $q=(x,y)$, by renumbering the bodies and choosing a proper coordinate system, we can always assume $x\in V_+$ and $y_i\neq y_j$ when $i\neq j$. Define
\begin{equation}\label{ma}
\text{Max}(q)=\max\{\frac{y_{i+1}-y_i}{x_{i+1}-x_i}| i=1,\dots,N-1\},
\end{equation}
\begin{equation}\label{mi}
\text{Min}(q)=\min\{\frac{y_{i+1}-y_i}{x_{i+1}-x_i}| i=1,\dots,N-1\},
\end{equation}
then $\text{Max}(q)\neq 0$ and $\text{Min}(q)\neq 0$.

The inequality in the following corollary allows us to check whether there exists a coordinate system in which the coordinates of the configuration can be arranged in the prescribed order.

\begin{cor} \label{MinMax}
    For a planar central configuration $q=(x,y)$ satisfying $x\in V_+$ under a proper coordinate system, $q$ is collinear if and only if 
    \begin{equation}\label{M1}
    -\text{Min}(q) \cdot \text{Max}(q)<1.
    \end{equation}
    
\end{cor}

\vspace{0.3cm}

Inequality (\ref{M1}) defines a neighborhood of the set of all collinear central configurations in the configuration space that contains no non-collinear central configurations.

\vspace{0.3cm}

\begin{proof}
    The inequality (\ref{M1}) holds for any $y\in V_+\cup V_-$. We just conclude with Proposition \ref{main}. 

    When $y\notin V_+\cup V_-$, we have $\text{Max}(q)>0$ and $\text{Min}(q)<0$. To prove the conclusion, we just need to prove that when inequality (\ref{M1}) holds, there is a rotation after which the coordinates of the configuration satisfy the prescribed ordering.
    
    We rotate the coordinate system $xOy$ by 
    $$
    \begin{pmatrix}
    \cos{\theta} & -\sin{\theta}\\
    \sin{\theta} & \cos{\theta}     
    \end{pmatrix},\quad \theta\in(0,\frac{\pi}{2}).
    $$
 Denote the new coordinates by $(x^*, y^*)$, then 
    $$
    x_i^*=x_i\cos{\theta}+y_i\sin{\theta},\quad y_i^*=-x_i\sin{\theta}+y_i\cos{\theta},\quad i=1,\dots,N.
    $$
We have
\begin{align*}
x^*\in V_+ & \iff x^*_i<x^*_{i+1}\\
           & \iff x_i\cos{\theta}+y_i\sin{\theta}<x_{i+1}\cos{\theta}+y_{i+1}\sin{\theta} \\
           &\iff \tan^{-1}{\theta}>-\frac{y_{i+1}-y_i}{x_{i+1}-x_i}\\
           &\iff \tan^{-1}{\theta}>\text{Max}\{-\frac{y_{i+1}-y_i}{x_{i+1}-x_i}|i=1,\dots,N-1\}\\
           &\iff \tan^{-1}{\theta}>-\text{Min}(q),
\end{align*}

\begin{align*}
    y^*\in V_+ &\iff y^*_i<y^*_{i+1}\\
    &\iff -x_i\sin{\theta}+y_i\cos{\theta}<-x_{i+1}\sin{\theta}+y_{i+1}\cos{\theta}\\
    &\iff \tan{\theta }<\frac{y_{i+1}-y_i}{x_{i+1}-x_i}\\
    &\iff \tan{\theta}<\text{Min}(q),
\end{align*}

\begin{align*}
    y^*\in V_- &\iff y^*_i>y^*_{i+1}\\
    &\iff -x_i\sin{\theta}+y_i\cos{\theta}>-x_{i+1}\sin{\theta}+y_{i+1}\cos{\theta}\\
    &\iff \tan{\theta }>\frac{y_{i+1}-y_i}{x_{i+1}-x_i}\\
    &\iff \tan{\theta}>\text{Max}(q).
\end{align*}
 Since $\tan{\theta}>0$, if $y^*\in V_+$, then $y\in V_+$, which is not the case under discussion. We denote by
$$
\Theta(q)=\{\theta\in (0,\frac{\pi}{2})| \tan^{-1}{\theta}>-\text{Min}(q);\ \tan{\theta}>\text{Max}(q)\}
$$
the admissible set of $\theta$ such that $x^*\in V_+$ and $y^*\in V_-$. Hence if $\Theta(q)\neq \emptyset$, the central configuration $q$ is collinear by Proposition \ref{main}. Since $\text{Max}(q)>0$ and $-\text{Min}(q)>0$, we have
\begin{itemize}
    \item $\Theta(x,y)=\emptyset$ when $-\text{Min}(q)>1$ and  $\text{Max}(q)>1$;
    \item $\Theta(x,y)\neq \emptyset$ when $0<-\text{Min}(q), \text{Max}(q)<1$;
    \item when $-\text{Min}(q)>1>\text{Max}(q)>0$ or $\text{Max}(q)>1>-\text{Min}(q)>0$, we have
    \begin{align*}
\Theta(q)\neq \emptyset &\iff 0<\arctan{(\text{Max}(q))}<\arccot{(-\text{Min}(q))}<\frac{\pi}{2}\\
&\iff 0<\text{Max}(q)<\tan{(\arccot{(-\text{Min}(q))})}=-\frac{1}{\text{Min}(q)}\\
&\iff 0<-\text{Min}(q)\cdot \text{Max}(q)<1.
\end{align*}
\end{itemize}
Since the second case is implied by $0<-\text{Min}(q)\cdot \text{Max}(q)<1$, we conclude.

\end{proof}

\section{Further Discussions}

Conley's Perpendicular Bisector Theorem, Moeckel's $45^{\circ}$ Theorem, Albouy and the first-named author's Disk Sector Theorem and our Theorem \ref{main} are all mass-independent criteria but are different in nature. On the one hand, for three bodies, we verified that the region defined by the inequality (\ref{M1}) in Corollary \ref{MinMax} is contained in the region defined by $\theta (q)$ in $45^{\circ}$ Theorem. While for more bodies, it is not clear. On the other hand, Theorem \ref{main} can be directly applied and generalized to higher dimensions. Some mass-dependent criteria are also established, for instance \cite{CC}.

After Theorem \ref{main}, a further question is whether there exist spatial central configurations which are `` almost coplanar ". In \cite{Moe}, Moeckel found a planar central configuration with 946 bodies that can be continued into spatial central configurations as masses vary. It has been proved this cannot happen for five bodies \cite{AS},\cite{AF}. Thus, it is reasonable to expect that criteria for ``almost coplanar" central configurations may depend on masses and the number of bodies.

Since a collinear central configuration also gives rise to a collinear balanced configuration with a given symmetric matrix, another further study is the generalization of Theorem \ref{main} to balanced configurations.

\vspace{0.3cm}

\noindent \textit{Acknowledgement.} We thank Alain Albouy for useful discussions. L. Z. is supported by fundamental research fund for the central universities and Liaoning revitalization programme.  

\noindent \textit{AI statement.} The proof is independent of AI. AI simplifies the computation of the equality (\ref{c}).




\end{document}